\documentclass[a4paper,draft]{amsproc}
\usepackage{amssymb}
\usepackage[hyphens]{url} 

\theoremstyle{plain}
 \newtheorem{theorem}{Theorem}[section]

 \newtheorem{corollary}{Corollary}[section]
\theoremstyle{definition}

\theoremstyle{remark}

 \numberwithin{equation}{section}

\renewcommand{\geq}{\geqslant}

\title[Linear processes/ M. Goubi]{On the the linear processes of a stationary time series $AR(2)$}



\author{\bfseries Mouloud  Goubi} 
\address{Mouloud Goubi\\
Department of Mathematics \\
University of UMMTO RP. 15000\\
Tizi-ouzou, Algeria\\
Laboratoire d'Alg\`ebre et Th\'eorie des Nombres, USTHB Alger}
\email{mouloud.goubi@ummto.dz}

\begin{document}

\vspace{18mm} \setcounter{page}{1} \thispagestyle{empty}

\begin{abstract}
Our aim in this work is to give explicit formula of the linear
processes solution of autoregressive time series $AR(2)$ with hint
of generating functions theory by using the Horadam numbers and
polynomials.
\end{abstract}

\maketitle

\section{Introduction}
The time series (see \cite{Bro02}) we observe are the realizations
of Random variables $Y_1,\cdots, Y_t$ which are a part of a larger
stochastic process $\left\{Y_t;\ t\in\mathbb{Z}\right\}$. Let
$\left\{Y_t;\ t\in\mathbb{Z}\right\}$ a zero-mean time series.
Consider $\mathcal{H}$ the Hilbert space spanned by the Random
variables $Y_t$ with inner product
\[\left\langle X,Y\right\rangle=E(X,Y)\]
and the norm
\[\|X\|=\sqrt{E(X^2)}.\]

The mean function of a time series is defined by $\mu(t)=E(Y_t)$ and
the auto-covariance function is given by
$\gamma(s,t)=cov\left(Y_s,Y_t\right)$. The mean and the
auto-covariance functions are fundamental parameters and it would be
useful to obtain sample estimates of them. For general time series
there are $2t +\frac{ t(t - 1)}{2}$ parameters associated with
$Y_1,\cdots, Y_t$ and it is not possible to estimate all these
parameters from $t$ data values. To make any progress we consider the times series is stationary.\\

On general a times series $Y_t$ is strictly stationary if for $k>0$
and any $t_1,\cdots,t_k$ of $\mathbb{Z}$, the distribution of
$\left(Y_{t_1},\cdots,Y_{t_k}\right)$ is the same as that for
$\left(Y_{t_1+u},\cdots,Y_{t_k+u}\right)$ for every $u$. If $Y_t$ is
stationary then $\mu(t)=\mu$ and $\gamma(s,t)=\gamma(s-t,0).$ We say
$Y_t$ is weakly stationary if $E(Y^2_t)<\infty,\ \mu(t)=0$ and
$\gamma(t+u,t)=\gamma(u,0).$ When time series are stationary it is
possible to simplify the parameterizations of the mean and
auto-covariance functions. In this case we can define the mean of
the series to be $\mu = E(Y_t)$ and the autocovariance function to
be $\gamma(u)=cov\left(Y_{t+u},Y_t\right)$.\\
If the random variables which make up $Y_t$ are uncorrelated, have
means $0$ and variance $\sigma^2$ (are so called white-noise
series). Then $Y_t$ is stationary with auto-covariance function
\begin{eqnarray}\label{chi}
\gamma(u)= \left\{
\begin{array}{ccc}
\sigma^2\ &\quad \textrm{ if }\ u=0, \\
0\  &\quad  \textrm{ otherwise}.
\end{array}
\right.
\end{eqnarray}
\section{Autoregressive series}

If the time series $Y_t$ satisfies the identity
\begin{equation}\label{eqTS}
Y_t=\phi_1Y_{t-1}+\cdots+\phi_pY_{t-p}+\epsilon_t
\end{equation}
where $\epsilon_t$ is white noise and $\phi_u$ are constants, then
$Y_t$ is called an autoregressive series of order $p$ and denoted by
$AR(p)$. These series are important, it explain how the next value
observed is a slight perturbation of a simple function of the most
recent observations. The solution; if exist in the form
\[Y_t=\displaystyle\sum_{u=0}^{\infty}\phi_u\epsilon_{t-u}\] is called a
linear processes with the condition
$\sum_{u=0}^{\infty}|\phi_n|^2<\infty$ to assure the convergence of
the last series in $\mathcal{H}$. The lag operator $L$ for a time
series $Y_t$ is defined by
\[L(Y_t)=Y_{t-1}\] and is linear. The last relation \eqref{eqTS} can
be written in the form
\begin{equation}\label{eqTS}
Y_t=\phi_1L(Y_{t})+\cdots+\phi_pL^p(Y_{t})+\epsilon_t.
\end{equation}
To found $Y_t$, we write
\[\big (1+\phi_1L+\cdots+\phi_pL^p\big )(Y_{t}) = \epsilon_t.\] If
the operator $1+\phi_1L+\cdots+\phi_pL^p$ is bijective we then write
\[Y_t=\left(1-\phi_1L-\cdots-\phi_pL^p\right)^{-1}\epsilon_t.\] The
linear processes of $AR(1)$ is completely known. We have
$Y_t=\left(1-L\right)^{-1}\epsilon_t$ then $Y_t=\sum_{u\geq0}\phi^{u
}_1\epsilon_{t-u}$ with the condition
$\sum_{u\geq0}|\phi_1|^{2n}<\infty$, which means that $|\phi_1|<1$.
An equivalent condition is that the root of the equation
$1-\phi_1z=0$ lies outside the unit circle in the complex
plane $\mathbb{C}$.\\

The $AR(2)$ model is defined by
$Y_t=\phi_1L(Y_t)+\phi_2L^2(Y_t)+\epsilon_t.$ Then we write
$Y_t=\left(1-\phi_1L-\phi_2L^2\right)^{-1}\epsilon_t.$ The
decomposition gives
\[1-\phi_1L-\phi_2L^2=\left(1-c_1L\right)\left(1-c_2L\right).\] We
can invert the operator if we can invert each factor separately.
This is possible if and only if $|c_1| < 1$ and $|c_2| < 1$, or
equivalently, if the roots of the polynomial $1-\phi_1z-\phi_2z^2$
lie outside the unit circle. But from the symmetric relations of the
roots we have $c_1+c_2=\phi_1$ and $c_1c_2=\phi_2$. The strong
conditions of convergence are $|\phi_1|<2$ and $|\phi_2|<1$.\\

In this work we are interested by times series belong to AR($2$)
respecting the conditions $|\phi_1|<2$ and $|\phi_2|<1$. The linear
processes $Y_t$ can be computed by two different methods. The first
expression; well-known by the statisticians which consist to compute
each $AR(1)$ separately and take the product. We reproduce it
directly by the following theorem.
\begin{theorem}\label{main1}
\begin{equation}\label{eqmain1}
Y_t=\sum_{u=0}^{\infty}\sum_{k=0}^{u}\phi^{u}_1\left(\phi_2/\phi_1\right)^k\epsilon_{t+k-u}\epsilon_{t-k}.
\end{equation}
\end{theorem}
The identity \eqref{eqmain1} is interesting, but the second member
of the equality contains the product of two white noise
$\epsilon_t$. Then it is not as the standard form of linear
processes. To escape this problem the second way is given by the
following theorem.
\begin{theorem}\label{main}
\begin{equation}\label{eqmain}
Y_t=\sum_{u=0}^{\infty}\left(\sum_{k=0}^{\left\lfloor\frac{u}{2}\right\rfloor}\binom{u-k}{k}\left(\phi_2/\phi^{2}_1\right)^k\right)\phi^{u}_1\epsilon_{t-u}.
\end{equation}
\end{theorem}
The condition $|\phi_1|<2$ and $|\phi_2|<1$ for the convergence of
the linear processes $Y_t$ on $\mathcal{H}$ can be replaced by the
condition
\[\sum_{u\geq0}\left(\sum_{k=0}^{\left\lfloor\frac{u}{2}\right\rfloor}\binom{u-k}{k}\left(\phi_2/\phi^{2}_1\right)^k\right)^2\phi^{2u}_1<\infty.\]
Letting $\phi_2=\phi^{2}_1$ in the expression \eqref{eqmain} Theorem
\ref{main}; the following corollary is immediate.
\begin{corollary}
\begin{equation}
Y_t=\sum_{u=0}^{\infty}\left(\sum_{k=0}^{\left\lfloor\frac{u}{2}\right\rfloor}\binom{u-k}{k}\right)\phi^{u}_1\epsilon_{t-u}.
\end{equation}
\end{corollary}
The symbol $\binom{u-k}{k}$ (see \cite{Com74}) occurs an important
place in combinatorics. It is a particular case of
$\binom{u-(k-1)l}{k}$ which is the number of k-blocks
$P\subset[u]=\left\{1,2,\cdots,n\right\}$ with the following
property; between two arbitrary points of $P$ are at least $l$
points of $[n]$ which do not belong to $P$.
\subsection{Proof of main results}
From the decomposition
$1-\phi_1L-\phi_2L^2=\left(1-c_1L\right)\left(1-c_2L\right).$ we
conclude that
\[Y_t=\left(\sum_{u\geq0}\phi^{u}_1\epsilon_{t-u}\right)\left(\sum_{u\geq0}\phi^{u}_2\epsilon_{t-u}\right)\]
Use Cauchy product of series to get the desired result
\eqref{eqmain1} Theorem \ref{main1}. For more details about this
technique we refer to \cite{Gou18}.\\

The proof of second theorem needs use techniques from generating
functions theory ( see \cite{Djor14}). We consider
\[\frac{1}{1-\phi_1xz-\phi_2z^2}=\sum_{u\geq0}A_u(x)t^u.\] where $A_u(x)$
 is a polynomial of degree $u$ to compute. Write
$f(x,z)=\frac{1}{1-\phi_1xz-\phi_2z^2}$ then
\[\left(1-\phi_1xz-\phi_2z^2\right)f(x,z)=1\]
and
\[\sum_{u\geq0}A_u(x)z^n-\phi_1x\sum_{u\geq0}A_n(x)z^{n+1}-\phi_2\sum_{u\geq0}A_u(x)z^{n+2}=1\]
Thus
\[\sum_{u\geq0}A_u(x)z^n-\phi_1x\sum_{u\geq1}A_{u-1}(x)z^{u}-\phi_2\sum_{u\geq2}A_{u-2}(x)z^{u}=1\]
and
\[A_0(x)+\left(A_1(x)-\phi_1xA_0(x)\right)z+\sum_{u\geq2}\left(A_u(x)-\phi_1xA_{u-1}(x)-\phi_2A_{u-2}(x)\right)z^u=1\]
Since this entire series is constant we deduce that $A_0(x)=1$,
$A_1(x)=\phi_1x$ and others are given by the recursion relation
\begin{equation}
A_u(x)=\phi_1xA_{u-1}(x)+\phi_2A_{u-2}(x).
\end{equation}
These polynomials are well-known and so called Horadam polynomials.
G. B. Djordjevi\'{c} and G. V. Milovanovi\'{c} (see \cite{Djor14})
provide the following explicit representation
\[A_u(x)=\sum_{k=0}^{\left\lfloor\frac{u}{2}\right\rfloor}\phi^{k}_2\frac{(u-k)!}{k!(u-2k)!}\left(\phi_1x\right)^{u-2k}.\]
Letting $x=1$ we obtain Horadam numbers $A_n$ these are written by
the form
\[A_u=\phi^{u}_1\sum_{k=0}^{\left\lfloor\frac{u}{2}\right\rfloor}\frac{(u-k)!}{k!(u-2k)!}\left(\phi_2/\phi^{2}_1\right)^{k}.\]
and admit for generating function
\begin{equation}\label{gf}
\frac{1}{1-\phi_1t-\phi_2t^2}=\sum_{u\geq0}A_ut^u.
\end{equation}
Numbers $A_n$ can be defined from the following recursion relation
\[A_u=\phi_1A_{u-1}+\phi_2A_{u-2},\ n\geq2\] with $A_0=1$ and
$A_1=\phi_1.$ Instead of $z$ we take the operator $L$ and it is
obvious to obtain the identity \eqref{eqmain} Theorem \ref{main}.


\begin{thebibliography}{99}

\bibitem{Bro02} P. J. Brockwell and R.A. Davis, {\it Introduction to Time Series and Forecasting}, Springer 2002.

\bibitem{Com74} L. Comtet, {\it Advanced combinatorics} D. Riedel Publishing Company Boston USA 1974.

\bibitem{Djor14} G. B. Djordjevic and G. V. Milovanovic, {\it Special Classes of  Polynomials}, University of Niš, Faculty of Technology, Leskovac, 2014.

\bibitem{Gou18} M. Goubi, {\it Successive derivatives of Fibonacci type polynomials of higher order in two variables}, Filomat {\bf 32}(4) (2018), pp. 5149--5159.

\end{thebibliography}
\end{document}